\documentclass[12pt]{article}

\usepackage{graphics,amsmath,amssymb,amsthm,mathrsfs,bm,color}
\usepackage{indentfirst}
\usepackage{amsmath}
\usepackage[titletoc,toc,title]{appendix}

\allowdisplaybreaks[2]
\newcommand{\lap}[1]{\Delta #1}

\newtheorem{theorem}{Theorem}[section]
\newtheorem{lemma}[theorem]{Lemma}

\newtheorem{definition}{Definition}[section]
\newtheorem{remark}{Remark}[section]

\numberwithin{equation}{section}

\begin{document}

\bibliographystyle{amsplain}

\title{A Hopf lemma and regularity for fractional $p$-Laplacians
\footnote{The first author is partially supported by the Simons Foundation Collaboration Grant for Mathematicians 245486. \newline
The second author is partially supported by NSFC-11571233. \newline
The third  author is supported by the
Fundamental Research Funds for the Central Universities (lzujbky-2017-it53).\newline
E-mail: wchen@yu.edu (W. Chen),~cli@colorado.edu (C. Li),~qishj15@lzu.edu.cn (S. Qi).}\\ }

\author{Wenxiong Chen,$^a$\quad Congming Li,$^b$\quad Shijie Qi$^c$\\
 \small{$^a$  Department of Mathematical Sciences, Yeshiva University }\\
\small{$^b$ School of Mathematics, Shanghai Jiao Tong University }\\
 \small{$^c$School of Mathematics and Statistics, Lanzhou University and Yeshiva University}}

\date{ }
\maketitle
\begin{abstract}
In this paper, we study qualitative properties of the fractional $p$-Laplacian. Specifically, we establish
a Hopf type lemma for positive weak super-solutions of the fractional $p-$Laplacian equation with Dirichlet condition. Moreover, an optimal condition is obtained to ensure $(-\triangle)_p^s u\in C^1(\mathbb{R}^n)$  for smooth functions $u$.
\vspace{0.3cm}\\
\small{\textbf{\emph{Keywords:}}\  Fractional $p-$Laplacian;  Hopf type lemma; regularity.}
\end{abstract}

\section{Introduction and main results}

The fractional $p-$Laplacian
is  defined by the singular integral
\begin{equation}\label{plaplacian}
\begin{aligned}
(-\triangle)_p^su(x)&:=&C_{n,s,p} P.V.\int_{\mathbb{R}^n}\frac{|u(x)-u(y)|^{p-2}(u(x)-u(y))}{|x-y|^{n+sp}}dy\\
&\equiv& C_{n,s,p}\lim_{\epsilon\to0}\int_{\mathbb{R}^n\setminus B_\epsilon(x)}\frac{|u(x)-u(y)|^{p-2}(u(x)-u(y))}{|x-y|^{n+sp}}dy,
\end{aligned}
\end{equation}
where $C_{n,s,p}$ is a positive constant depending only on $n,$ $s$, and $p$,  $s\in (0,1)$, and $p>1$. Denote
$$L_{sp}(\mathbb{R}^n):=\left\{u\in L_{\text{loc}}^1(\mathbb{R}^n)\Big|\int_{\mathbb{R}^n}\frac{|u(x)|^{p-1}}{(1+|x|)^{n+sp}}dx<\infty\right\}.$$
If $u\in C_{\text{loc}}^{1,1}\cap L_{sp}(\mathbb{R}^n)$, then \eqref{plaplacian} is well defined.
Clearly, when $p=2$, \eqref{plaplacian} becomes the fractional
Laplacian which arises in many fields such as phase transitions, flame propagation, stratified materials and others (see \cite{ab,bgr,sirevaldinoci}). In particular, the fractional Laplacian can be understood as the infinitesimal generator of a stable Levy process (see \cite{valdinoci}). The fractional $p-$Laplacian also has many applications, for instance, it is used to study the non-local ``Tug-of-War" game (see \cite{BCF1,BCF2,ishiinakamura}). The interest on these  nonlocal operators continues to grow  in recent years. We refer to \cite{ms} for the recent progress on these nonlocal operators.

Due to the non-locality of these kinds of operators, many traditional methods in studying the local differential operators no longer work.
To overcome this difficulty,
Cafarelli and Silvestre \cite{cs} introduced the {\em extension method}
which turns nonlocal problems involving the fractional Laplacian  ($p=2$) into local ones in higher dimensions, then the classical theories for local elliptic partial differential equations can be applied.
We refer to \cite{bcps,chenzhu} and references therein for broad applications of this method.

Another useful method to study the fractional Laplacian is the {\em integral equations method}, which turns a given fractional Laplacian equation into its equivalent integral equation, and then various properties of the original equation can be obtained by investigating the integral equation, see \cite{cfy,chenlima,zccy} and references therein.

However, so far as we know, there has neither been any  extension method nor the integral equations method that work for the fractional $p-$Laplacian equation  when $p\neq2$. The nonlinearity, the singularity $(1<p<2)$ and degeneracy $(p>2)$ of the operator $(-\triangle)_p^s$ render many powerful methods to study the fractional Laplacian $(p=2)$ no longer effective.

Recently, Chen et al. have developed a {\em direct method of moving planes} to investigate the nonlocal problems, which can be used to study not only the fractional Laplacian but also the fully nonlinear nonlocal operator
$$F_\alpha(u(x)):=C_{n,\alpha}\lim_{\epsilon\to0}\int_{\mathbb{R}^n\setminus B_\epsilon(x)}\frac{G(u(x)-u(y))}{|x-y|^{n+\alpha}}dy,$$
where $\alpha>0$, $G:$ $\mathbb{R}\to \mathbb{R}$ is a locally Lipschitz continuous function. The fractional $p$-Laplacian is a special case in  which $G(t)=|t|^{p-2}t$ and $\alpha = s p$.
This direct method has been successfully applied to obtain symmetry, monotonicity,  nonexistence and other qualitative properties of solutions for various nonlocal problems, see e.g., \cite{chenli1,chenlizhang,chenlili,chenliliyan1,chenliliyan2}.

In the present paper, we will continue to study qualitative properties for fractional $p$-Laplacian.
We will establish a Hopf type lemma in general domains for super solutions to fractional $p$-Laplacian equations with a Dirichlet condition; and for any given smooth function $u$, we will obtain an optimal condition for $(-\Delta)^s_p u$ to be continuously differentiable.

It is well-known that the Hopf lemma is a very powerful tool in the study of various differential equations. For example, it has been successfully used in the ``second" step of the moving planes method.

In the case of fractional Laplacian ($p=2$), Fall and Jarohs \cite[Proposition 3.3 ]{fall-jarohs}  proved a Hopf lemma for the  entire antisymmetric supersolution of the problem
\begin{equation}\label{firsthopf}
(-\triangle)^su(x)= c(x)u(x)\quad \text{in}~ \Omega.
\end{equation}
 Greco and Servadei \cite{grecoservadei} obtained a Hopf type lemma to \eqref{firsthopf}
under the assumptions that $c(x)\leq0$ and $\Omega\subset \mathbb{R}^n$ is a bounded domain. Chen and Li \cite{chenli2} established a Hopf lemma for anti-symmetric function on a half space through a rather delicate analysis. More recently,
Jin and Li \cite{jinli} extended the results in \cite{chenli2} to the fractional $p-$Laplacian with $p>3$
for positive anti-symmetric functions on the boundary of  a half space.  In this paper, we shall establish a Hopf type lemma  for the positive weak supersolution of  \eqref{PLD} on the boundary of more general domains.

Before stating our main results, we first introduce some definitions on fractional Sobolev spaces, and one can see \cite{dpv,iannizzotto} for more details. For any domain $\Omega\subset\mathbb{R}^n$ with smooth boundary, define
$$W^{s,p}(\Omega):=\left\{u\in L^p(\Omega)\Big| \int_\Omega\int_\Omega\frac{|u(x)-u(y)|^p}{|x-y|^{n+sp}}dxdy<\infty\right\}$$
equipped with the norm
$$||u||_{W^{s,p}(\Omega)}:=||u||_{L^p(\Omega)}+\left(\int_\Omega\int_\Omega\frac{|u(x)-u(y)|^p}{|x-y|^{n+sp}}dxdy\right)
^{\frac{1}{p}},$$
and
$$W^{s,p}_0(\Omega):=\left\{u\in W^{s,p}(\mathbb{R}^n)\mid u=0 \text{~in~}\mathbb{R}^n\backslash\Omega\right\}.$$
If $\Omega\subset \mathbb{R}^n$ is bounded, set
$$\widetilde{W}^{s,p}(\Omega):=\left\{u\in L_{loc}^p(\mathbb{R}^n)\Big| \; \exists U\supset\supset\Omega
\text{ such that }||u||_{W^{s,p}(U)}+\int_{\mathbb{R}^n}\frac{|u(x)|^{p-1}}{(1+|x|)^{n+sp}}dx<\infty\right\}.$$
If $\Omega\subset \mathbb{R}^n$ is unbounded, set
$$\widetilde{W}^{s,p}_{\text{loc}}(\Omega):=\{u\in L^p_{\text{loc}}(\mathbb{R}^n) \mid  u\in \widetilde{W}^{s,p}(\Omega^\prime)\text{~for any~}
\Omega^\prime\subset\subset\Omega\}.$$
Next, we present two definitions of solutions to fractional $p-$Laplacian equation with Dirichlet condition
\begin{equation}\label{PLD}
\begin{cases}
(-\triangle)_p^su(x)=f(x)&\text{~in~}\Omega,\\
u=0 &\text{~in~}\mathbb{R}^n\backslash\Omega.
\end{cases}
\end{equation}

\begin{definition}
We say that $u$ is a classical supersolution (subsolution) of the Dirichlet problem \eqref{PLD}, if \eqref{plaplacian} is well-defined for any $x\in\Omega$,  moreover, there hold
\begin{equation}
\begin{cases}
(-\triangle)_p^su(x)\geq (\leq) f(x)&\text{~in~}\Omega,\\
u\geq (\leq)0 &\text{~in~}\mathbb{R}^n\backslash\Omega.
\end{cases}
\end{equation}
Furthermore, if $u$ is both a  supersolution and a  subsolution of \eqref{PLD}, then we say it is a  solution to \eqref{PLD}.
\end{definition}

\begin{definition}
Let $f\in W^{-s,p^\prime}(\Omega),$ we say that $u\in \widetilde{W}^{s,p}(\Omega)$  is a  weak supersolution of \eqref{PLD}, if there hold
$$(u+\epsilon)^-\in W^{s,p}_0(\Omega)\text{~for any~}\epsilon>0,$$
and
$$C_{n,s,p}\int_{\mathbb{R}^n}\int_{\mathbb{R}^n}\frac{|u(x)-u(y)|^{p-2}(u(x)-u(y))(\phi(x)-\phi(y))}{|x-y|^{n+sp}}dxdy\geq
\langle f,\phi\rangle$$
for any $\phi\in W^{s,p}_0(\Omega)$ with $\phi\geq0$ in $\Omega.$
The weak subsolution can be defined similarly.
 Moreover, if $u$ is both a weak supersolution and a weak subsolution of \eqref{PLD}, then we say it is a weak solution to \eqref{PLD}.
\end{definition}

One of our main results is

\begin{theorem}\label{hopf}
 Let $\Omega\subset\mathbb{R}^n $ be a domain with $C^{1,1}$ boundary. If it is bounded, we assume $u\in \widetilde{W}^{s,p}(\Omega)\cap C(\Omega)$; if it is unbounded, we assume
 $u\in \widetilde{W}_{\text{loc}}^{s,p}(\Omega)\cap C(\Omega)$.

 Suppose

\begin{equation}\label{Q}
\begin{cases}
(-\triangle)_p^su\geq 0 \quad &{\rm\text{in}}~\Omega,\\
u>0 &\text{in}~\Omega,\\
u=0 &\text{in}~\mathbb{R}^n\backslash\Omega
\end{cases}
\end{equation}
in the weak sense, then
$$\mathop{\liminf}\limits_{d(x)\rightarrow0}\frac{u(x)}{d^s(x)}>0,$$
where $d(x):=dist(x,\Omega^c)$.
\end{theorem}

The other main result is concerning the regularity of $(-\Delta)^s_p u$.

The regularity of solutions of the fractional $p-$Laplacian equations has attracted considerable attention in recent years, and it has been well understood for the fractional Laplacian equations ($p=2$).
Specifically, the Schauder interior estimate of the solution is similar to that of the Poisson equation (associated with the regular Laplacian), which states roughly that if $f \in C^{\gamma}(\Omega)$ and $u\in C_{\text{loc}}^{1,1}(\Omega)\cap L_{2s}(\mathbb{R}^n)$ is a solution of
\begin{equation}\label{FD}
\begin{cases}
(-\triangle)^su(x)=f(x)&\text{~in~}\Omega,\\
u=0 &\text{~in~}\mathbb{R}^n\backslash\Omega,
\end{cases}
\end{equation}
then the regularity of the solution $u$ can be raised by the order of ${2s}$ in any proper subset of $\Omega$, the same order as the operator $(-\lap)^{s}$. By introducing the proper weighted H$\ddot{o}$lder norms as in the case of Poisson equations, one shall be able to control a weighted $C^{{2s} +\gamma}$ norm of $u$ in $\Omega$ in terms of another weighted $C^{\gamma}$ norm of $f$ in $\Omega$. However, when considering the  regularity of the solution up to the boundary, the situation in the fractional order equation is quite different from that in the integer order equation (when $s=1$, the Poisson equation). In fact, Ros-Oton and Serra \cite{otonserra} proved that if
$u\in C_{\text{loc}}^{1,1}(\Omega)\cap L_{2s}(\mathbb{R}^n)$ is a solution of \eqref{FD} with $f\in L^\infty(\Omega),$
then $u$ is $C^s$ up to the boundary; and this is optimal in general.  Later,
Chen et al. \cite{chenlima} proved the similar results by a simpler method.

For the fractional $p-$Laplacian, the study of the regularity becomes quite complicated.  So far there are very few results.
Di Castro and Kuusi \cite{ck} showed that if $u\in\widetilde{W}^{s,p}(\Omega)$ satisfies $(-\triangle)^s_pu=0$ in $\Omega,$ then $u$ is locally $\gamma$-H$\ddot{o}$lder continuous for  small $\gamma$. Brasco et al. \cite{brasco} established a higher H$\ddot{o}$lder regularity for the fractional $p-$Laplacian equation in the superquadratic case $(p>2)$. Indeed, the authors have verified that if $u\in W_{\text{loc}}^{s,p}(\Omega)\cap L_{sp}(\mathbb{R}^n)$ is a local weak solution of
 \begin{equation}\label{pe}
 (-\triangle)_p^su=f\quad\text{in}~\Omega,
 \end{equation}
 where $f\in L^q_{\text{loc}}(\Omega)$ with
 \begin{equation}
 \begin{cases}
 q>\frac{n}{sp}\quad &\text{if}~sp\leq n,\\
 q\geq1 &\text{if}~sp>n,
 \end{cases}
 \end{equation}
then $u\in C^\delta_{\text{loc}}(\Omega)$ for every $0<\delta<\Theta(n,s,p,q)$ with
$$\Theta(n,s,p,q)=\min\left\{\frac{1}{p-1}\left(sp-\frac{n}{q}\right),1\right\}.$$
Iannizzotto et al. \cite{iannizzotto} proved that the solutions of \eqref{pe} with $f\in L^\infty(\Omega)$ belong to
$C^\alpha(\overline{\Omega})$ for some $\alpha\in (0,s].$

Concerning the regularities of $(-\triangle)_p^su$ for a given smooth function $u$, there are more substantial technical difficulties than the local case.

For the fractional Laplacian $(-\triangle)^s,$ Silvestre \cite{silvestre} has made a comprehensive investigation. More specifically, he has verified that
if $u\in L_{2s}(\mathbb{R}^n)\cap C^{2s+\epsilon}$ (or $C^{1,2s+\epsilon-1}$ if $s>1/2$) for some $\epsilon>0$ in an open set $\Omega,$ then $(-\triangle)^su$ is a continuous function in $\Omega$ for $s\in (0,1).$ Furthermore, if $u\in C^{k,\alpha}$ and  $k+\alpha-2s$ is not an integer, then $(-\triangle)^su\in C^{l,\beta}$, where $l$ is the integer part of $k+\alpha-2s$ and $\beta=k+\alpha-2s-l.$

While for the fractional $p-$Laplacian, the singularity ($0<p<2$) and degeneracy ($p>2$) of operator $(-\triangle)_p^s$ make it more complex.

For example, even for the local operator $\triangle_p$ and the sufficient smooth function $u(x)=x^2$ in $\mathbb{R},$ $-\triangle_p u(0)=\infty$ if $1<p<2,$ and  $(-\triangle_pu)^\prime(0)=\infty$ if $1<p<3$ and $p\neq2.$

In this paper, we shall consider the differentiability of $(-\triangle)_p^su$ for $p>2$ and establish an optimal condition such that $(-\triangle)_p^su\in C^1(\mathbb{R}^n)$. Specifically, we prove that

\begin{theorem}\label{differeniability}
Let $p>2,$   $u\in C_{\text{loc}}^3(\mathbb{R}^n)\cap L_{sp}(\mathbb{R}^n)$ and $|\nabla u|\in L_{sp}(\mathbb{R}^n).$ If $p>\frac{3}{2-s},$
then $(-\triangle )_p^su\in C^1(\mathbb{R}^n).$
\end{theorem}

The condition $p>\frac{3}{2-s}$ is optimal as shown in the following

\begin{theorem}\label{opposite}
Let $p>2,$ $u(x)=\eta(x)x^2$ in $\mathbb{R},$ where $\eta\in C_0^\infty(-2,2)$ is even and satisfies $0\leq \eta(x)\leq 1$ in $\mathbb{R},$ $\eta(x)=1$ in $(-1,1)$ and $|\eta^\prime(x)|\leq 1$ in $\mathbb{R}.$ If $p<\frac{3}{2-s},$
then
\begin{equation}\label{example}
\lim_{x\to0+}\left| \left((-\triangle)_p^su\right)^\prime(x)\right|=\infty.
\end{equation}
And if $p=\frac{3}{2-s},$ then
\begin{equation}\label{examplecase2}
\lim_{x\to 0+}\left((-\triangle)^s_pu\right)^\prime(x)\neq\left((-\triangle)^s_pu\right)^\prime(0).
\end{equation}
\end{theorem}

The rest of the paper is organized as follows. Section 2 is devoted to establishing the Hopf type lemma for the positive solution
of $\eqref{Q}$. In section 3, we first prove the differentiability of $(-\triangle)_p^su$ under the condition $p>\frac{3}{2-s}$. Then we show that this condition is optimal by giving a counterexample when $p\leq\frac{3}{2-s}$.
In { the Appendix, we state some results in \cite{iannizzotto} used in the present paper for convenience.}

\section{Hopf type lemma}

In this section, we prove the Hopf type lemma for the positive weak solution of \eqref{Q} by constructing a suitable subsolution.

\begin{proof}[\bf Proof of Theorem \ref{hopf}]
For any given $x_0\in \partial\Omega,$ it follows from the $C^{1,1}$ property of the boundary of $\Omega$ that there exist $x_1\in\Omega$ on the normal line to $\partial \Omega$ at $x_0$ and a positive constant $\alpha$ such that $B_\alpha(x_1)\subset \Omega,$ $\overline{B}_\alpha(x_1)\cap\partial\Omega=x_0 $ and dist$(x_1,\Omega^c)=|x_0-x_1|.$ Without loss of the generality, we suppose that $x_0$ is the origin, $\alpha=1$ and $x_1=e_n$ with $e_n=(0,\cdots,1)$ the last vector of the canonical basis of $\mathbb{R}^n.$
Let $r\in\left(0,\frac{1}{3\sqrt{5}}\right)$ be a constant, $O$ denote the origin and   $\eta\in C^2(\mathbb{R}^n)$ satisfy that
\begin{equation}\label{eta}
\eta(X)=1~\text{in}~ B_{2r}(O),~~\eta(X)=0~\text{in}~ B_{3r}^c(O), ~~\text{and}~ |\nabla \eta|\leq\frac{1}{r}~\text{in}~\mathbb{R}^n.
\end{equation}
Now, define $\Psi: \mathbb{R}^n\to\mathbb{R}^n$ as
\begin{equation}\label{pa}
\Psi(X)=X+\left(1-X_n-\sqrt{((1-X_n)^2-|X^\prime|^2)_+}\right)\eta(X)e_n\quad \text{for any}~ X\in \mathbb{R}^n,
\end{equation}
where $X=(X^\prime,X_n).$  Clearly, it follows from
\eqref{eta} and \eqref{pa} that $\Psi(X)=X$ for any $X\in B_{3r}^c(O).$
Since $r<\frac{1}{3\sqrt{5}},$ we have
\begin{equation}\label{X}
1-X_n\geq 2|X^\prime|\quad \text{for any}~ X\in B_{3r}(O),
\end{equation}
which implies
\begin{equation}\label{pa1}
 \Psi(X)=
 \begin{cases}
 X+\left(1-X_n-\sqrt{(1-X_n)^2-|X^\prime|^2}\right)\eta(X)e_n~&\text{for any}~ X\in B_{3r}(O),\\
 X&\text{for any}~ X\in  B_{3r}^c(O).
\end{cases}
\end{equation}
{Next we show the following two claims.}

{\bf Claim 1.} $\Psi$ is a $C^{1,1}$ diffeomorphism of $\mathbb{R}^n.$  We firstly show that $\Psi$ is a bijection in $B_{3r}(O).$
Noting that if there exist $X=(X^\prime, X_n),$  $Y=(Y^\prime, Y_n)\in \mathbb{R}^n$ such that  $\Psi(X)=\Psi(Y)$, then $X^\prime=Y^\prime.$
For any given $X^\prime\in\mathbb{R}^{n-1}$ with $|X^\prime|\leq3r,$ define
$$h:\{X_n\in\mathbb{R}|(X^\prime,X_n)\in B_{3r}(O)\}\to\mathbb{R}$$ as
\begin{equation}\label{h}
h(X_n)=X_n+\left(1-X_n-\sqrt{(1-X_n)^2-|X^\prime|^2}\right)\eta(X).
\end{equation}
Now we  show that $h$ is strictly monotone. Direct calculation implies that
\begin{equation}\label{hprime}
h^\prime(X_n)=1+\frac{|X^\prime|^2}{1-X_n+\sqrt{(1-X_n)^2-|X^\prime|^2}}\frac{\partial\eta}{\partial X_n}(X)
+\left(-1+\frac{1-X_n}{\sqrt{(1-X_n)^2-|X^\prime|^2}}\right)\eta(X).
\end{equation}
Noting that if $|X^\prime|=0$ then $h^\prime=1.$  Furthermore, if $|X^\prime|\neq0,$ the last term in \eqref{hprime} is positive and the second term can be rewritten as
\begin{equation*}
\frac{|X^\prime|}{1-X_n+\sqrt{(1-X_n)^2-|X^\prime|^2}}\cdot\frac{\partial \eta}{\partial X_n}(X)|X^\prime|:=J_1\cdot J_2.
\end{equation*}
It then follows from $|X^\prime|\leq 3r$ and \eqref{eta} that $|J_2|\leq 3.$ Moreover, thanks to \eqref{X}, we can verify that
$J_1\leq\frac{1}{2+\sqrt{3}}.$ Consequently, there holds
$$h^\prime(X_n)\geq 1-\frac{3}{2+\sqrt{3}}>0,$$
that is, $h$ is strictly increasing in $\{X_n\in\mathbb{R}|(X^\prime,X_n)\in B_{3r}(O)\}$.  This together with \eqref{pa1} shows that $\Psi\in C^{1,1}(\mathbb{R}^n,\mathbb{R}^n)$ is a diffeomorphism of $\mathbb{R}^n.$

{\bf Claim 2.} There holds
\begin{equation}\label{papro}
B_1(e_n)\cap B_r(0)\subset\Psi(B_{\sqrt{2}r}(O)\cap \mathbb{R}^n_+).
\end{equation}
{Indeed,} it suffices to show that  for any $x\in B_1(e_n)\cap B_r(0),$ there is $X\in B_{\sqrt{2}r}(O)\cap \mathbb{R}^n_+$ such that  $\Psi(X)=x.$
To this end, choose
$$X=\left(x^\prime, 1-\sqrt{(1-x_n)^2+|x^\prime|^2}\right).$$
Then
\begin{align*}
|X|^2&=|x^\prime|^2+1-2\sqrt{(1-x_n)^2+|x^\prime|^2}+(1-x_n)^2+|x^\prime|^2\\
&\leq 2|x^\prime|^2+1-2(1-x_n)+(1-x_n)^2\\
&\leq 2|x|^2\\
&\leq 2r^2,
\end{align*}
which shows that $X\in B_{\sqrt{2}r}(O)\cap \mathbb{R}^n_+$. It then follows from \eqref{pa1} that
\begin{align*}
\Psi(X)&=X+\left(1-X_n-\sqrt{(1-X_n)^2-|X^\prime|^2}\right)e_n\\
&=\left(x^\prime, 1-\sqrt{(1-x_n)^2+|x^\prime|^2}\right)+\left(\sqrt{(1-x_n)^2+|x^\prime|^2}-(1-x_n)\right)e_n\\
&=x.
\end{align*}
The claim is then verified.

Now, we define $\rho:\mathbb{R}^n\to \mathbb{R}$ as
\begin{equation}\label{rho}
\rho(x)={\rm\text{dist}}(x,B^c_1(e_n))\quad {\rm\text{in}}~\mathbb{R}^n.
\end{equation}
Then there holds
\begin{equation}\label{paprop}
\rho(\Psi(X))=(X_n)_+ \quad \text{for any}~X\in B_{2r}(O).
\end{equation}
Indeed, It follows from \eqref{eta} and \eqref{pa1} that
$$\Psi(X)=\left(X^\prime,1-\sqrt{(1-X_n)^2-|X^\prime|^2}\right)\quad \text{for any}~X\in B_{2r}(O).$$
By a direct calculation, we see that $\Psi(X)\in \partial B_{1-X_n}(e_n)$ for any $X\in B_{2r}(O),$ which
together with \eqref{rho} leads to \eqref{paprop}. Next we show by a similar calculation as in \cite{iannizzotto} that there exists a positive constant $C_1$ such that
\begin{equation}\label{rhos}
(-\triangle)_p^s\rho^s(x)\leq C_1\quad \text{in}~B_1(e_n)\cap B_r(0).
\end{equation}
Indeed, thanks to lemma \ref{generaldef}, we only need to show that there exists $f\in L^\infty(B_1(e_n)\cap B_r(0))$ such that
$$\mathop{\lim}\limits_{\epsilon\to0}\int_{\{|\Psi^{-1}(x)-\Psi^{-1}(y)|>\epsilon\}}\frac{G(\rho^s(x)-\rho^s(y))}{|x-y|^{n+ps}}dy
=f \quad {in}~L^1(B_1(e_n)\cap B_r(0)),$$
 where $G(t)=|t|^{p-2}t$ for any $t\in \mathbb{R}.$
Making a change of variables $X=\Psi^{-1}(x),$ then for any $x\in B_1(e_n)\cap B_r(0),$ there  exists
$X\in B_{\sqrt{2}r}(O)\cap\mathbb{R}^n_+$ such that $\Psi(X)=x$  and

\begin{align}\label{hopfj12}
&\int_{\{|\Psi^{-1}(x)-\Psi^{-1}(y)|>\epsilon\}}\frac{G(\rho^s(x)-\rho^s(y))}{|x-y|^{n+ps}}dy\notag\\
=&\int_{B_\epsilon^c(X)}\frac{G\left(\rho^s(\Psi(X))-\rho^s(\Psi(Y))\right)}{|\Psi(X)-\Psi(Y)|^{n+ps}}J_{\Psi}(Y)dY\notag\\
=&\int_{B_\epsilon^c(X)\cap B_{2r}(O)}\frac{G\left((X_n)_+^s-(Y_n)_+^s\right)}{|\Psi(X)-\Psi(Y)|^{n+ps}}J_{\Psi}(Y)dY
+\int_{B_{2r}^c(O)}\frac{G\left(\rho^s(\Psi(X))-\rho^s(\Psi(Y))\right)}{|\Psi(X)-\Psi(Y)|^{n+ps}}J_{\Psi}(Y)dY\notag\\
=&\int_{B_\epsilon^c(X)}\frac{G\left((X_n)_+^s-(Y_n)_+^s\right)}{|\Psi(X)-\Psi(Y)|^{n+ps}}J_{\Psi}(Y)dY\notag\\
&+\int_{B_{2r}^c(O)}\frac{G\left(\rho^s(\Psi(X))-\rho^s(\Psi(Y))\right)-G\left((X_n)_+^s-(Y_n)_+^s\right)}
{|\Psi(X)-\Psi(Y)|^{n+ps}}J_{\Psi}(Y)dY\notag\\
:=&J_1(X)+J_2(X),
\end{align}
where the second equality follows from \eqref{papro} and \eqref{paprop}. Noting that $\Psi$ is a $C^{1,1}$ diffeomorphism of $\mathbb{R}^n$ and $\Psi=I$ in $B_{3r}^c(O),$ Lemma \ref{changeofvariable}  then yields that
there exists $f_1\in L^\infty(B_1(e_n)\cap B_r(0))$ such that
\begin{equation}\label{j1}
\begin{aligned}
\mathop{\lim}\limits_{\epsilon\to0}J_1(X)&=
\lim_{\epsilon\to0}\int_{B_\epsilon^c(X)}\frac{G\left((X_n)_+^s-(Y_n)_+^s\right)}{|\Psi(X)-\Psi(Y)|^{n+ps}}J_{\Psi}(Y)dY\\
&=\lim_{\epsilon\to0}\int_{\{|\Psi^{-1}(x)-\Psi^{-1}(y)|>\epsilon\}}
\frac{G\left((\Psi^{-1}(x)\cdot e_n)_+^s-(\Psi^{-1}(y)\cdot e_n)_+^s\right)}{|x-y|^{n+sp}}dy\\
&=f_1(\Psi(X)) \quad\text{in} ~L^1\left(\Psi^{-1}(B_1(e_n)\cap B_r(0))\right).
\end{aligned}
\end{equation}
Thanks to \eqref{papro},  there exists a positive constant $C$ such that for any $x\in B_1(e_n)\cap B_r(0)$ and $ Y\in B_{2r}^c(O),$
$$|X-Y|\geq C(1+|Y|) .$$
Hence, we have
\begin{equation*}
\begin{aligned}
|J_2(X)|&\leq\int_{B_{2r}^c(O)}\frac{\left|\rho^s(\Psi(X))-\rho^s(\Psi(Y))\right|^{p-1}+\left|(X_n)_+^s-(Y_n)_+^s\right|^{p-1}}
{|\Psi(X)-\Psi(Y)|^{n+ps}}J_{\Psi}(Y)dY\\
&\leq C_\Psi\int_{B_{2r}^c(O)}\frac{1}{|X-Y|^{n+s}}dY\\
&\leq C_\Psi\int_{B_{2r}^c(O)}\frac{1}{(1+|Y|)^{n+s}}dY\\
&\leq C_\Psi,
\end{aligned}
\end{equation*}
where the notation $C_\Psi$ above may denote different positive constants. This together with \eqref{hopfj12}
and \eqref{j1} shows that
$$\mathop{\lim}\limits_{\epsilon\to0}
\int_{B_\epsilon^c(X)}\frac{G\left(\rho^s(\Psi(X))-\rho^s(\Psi(Y))\right)}{|\Psi(X)-\Psi(Y)|^{n+ps}}J_{\Psi}(Y)dY=f_1(\Psi(X))+ J_2(X)$$
in $L^1\left(\Psi^{-1}(B_1(e_n)\cap B_r(0))\right)$, with $f_1\circ\Psi$ and $J_2$ belong to $L^\infty\left(\Psi^{-1}(B_1(e_n)\cap B_r(0))\right).$ It then follows that
$$\mathop{\lim}\limits_{\epsilon\to0}\int_{\{|\Psi^{-1}(x)-\Psi^{-1}(y)|>\epsilon\}}\frac{G(\rho^s(x)-\rho^s(y))}{|x-y|^{n+ps}}dy
=f_1(x) +J_2\circ\Psi^{-1}(x) \quad \text{in}~L^1(B_1(e_n)\cap B_r(0)).$$
Consequently, \eqref{rhos} follows.

Now let $D\subset\subset B_1^c(e_n)\cap\Omega$ be a bounded smooth domain, and $\beta>0$
be a positive constant to be determined below. Set
\begin{equation}\label{sub}
\underline{u}(x)=\beta\rho^s(x)+\chi_D(x)u(x),
\end{equation}
where $\rho$ is defined by \eqref{rho}, and $\chi_D$ is the  characteristic function of $D,$  namely,
\begin{equation*}
\chi_{D}(x)=
\begin{cases}
1,\quad &x\in D,\\
0, &x\notin D.
\end{cases}
\end{equation*}
It follows from $D\subset\subset B_1^c(e_n)\cap\Omega$ that there is a positive constant $C_D$ such that
\begin{equation}\label{cd}
|x-y|\geq C_D\quad \text{for any}~ x\in B_1(e_n), ~y\in D.
\end{equation}
For any $x\in B_1(e_n)\cap B_r(0),$ direct calculation (we omit the term `$C_{n,s,p}\mathop{\lim}\limits_{\epsilon\to 0}$' in the following calculation for convenience) shows that
\begin{align}\label{subsolution}
(-\triangle)_p^s\underline{u}(x)
&=\int_{B^c_\epsilon(x)}\frac{G(\underline{u}(x)-\underline{u}(y))}{|x-y|^{n+ps}}dy\notag\\
&=\int_{B^c_\epsilon(x)}\frac{G\left(\beta\rho^s(x)-\beta\rho^s(y)-\chi_D(y)u(y)\right)}{|x-y|^{n+ps}}dy\notag\\
&=\int_{B^c_\epsilon(x)\cap B_1(e_n)}\frac{G\left(\beta\rho^s(x)-\beta\rho^s(y)\right)}{|x-y|^{n+ps}}dy
+\int_{B^c_{1}(e_n)}\frac{G\left(\beta\rho^s(x)-\chi_D(y)u(y)\right)}{|x-y|^{n+ps}}dy\notag\\
&=\int_{B^c_\epsilon(x)}\frac{G\left(\beta\rho^s(x)-\beta\rho^s(y)\right)}{|x-y|^{n+ps}}dy\notag\\
&~~+\int_{B^c_{1}(e_n)}\frac{G\left(\beta\rho^s(x)-\chi_D(y)u(y)\right)-G\left(\beta\rho^s(x)-\beta\rho^s(y)\right)}
{|x-y|^{n+ps}}dy\notag\\
&=\beta^{p-1}(-\triangle)_p^s\rho^s(x) +\int_{D}\frac{G\left(\beta\rho^s(x)-u(y)\right)-G\left(\beta\rho^s(x)\right)}{|x-y|^{n+ps}}dy\notag\\
&\leq \beta^{p-1}C_1 +\overline{C}_DA(x),
\end{align}
where
\begin{equation}\label{A}
A(x)=\int_DG\left(\beta\rho^s(x)-u(y)\right)-G\left(\beta\rho^s(x)\right)dy,
\end{equation}
and the last inequality holds due to \eqref{rhos}. Let
$$M_0=\mathop{\min}\limits_{x\in D}u(x)>0,\quad \beta\leq \frac{1}{2}M_0,$$ then it follows from the monotonicity of $G$ that
\begin{equation}\label{A1}
A(x)\leq\int_{D}G(\beta\rho^s(x)-u(y))dy\leq\int_{D}G\left(\frac{1}{2}M_0-M_0\right)dy=-\left(\frac{1}{2}\right)^{p-1}M_0^{p-1}|D|.
\end{equation}
It then follows from \eqref{subsolution}, \eqref{A} and \eqref{A1} that
\begin{equation*}
(-\triangle)_p^s\underline{u}(x)\leq M_1\beta^{p-1}-M_2\quad \text{in}~B_1(e_n)\cap B_r(0),
\end{equation*}
where $M_1 $ and $M_2$ are some positive constants. In view of \eqref{rho} and \eqref{sub}, there holds
$$\underline{u}(x)\leq u(x)\text{~~in~~}B_1^c(e_n).$$
Let
$$M_3=\mathop{\inf}\limits_{x\in B_1(e_n)\cap B_r^c(0)}u(x)>0,\quad \beta<\frac{1}{2}\min\left\{M_0,M_3,\left(\frac{M_2}{M_1}\right)^{\frac{1}{p-1}}\right\},$$
then we have
\begin{equation*}
\begin{cases}
(-\triangle)_p^s\underline{u}< 0\quad& \text{in}~B_1(e_n)\cap B_r(0),\\
\underline{u}(x)\leq u(x) &\text{in}~(B_1(e_n)\cap B_r(0))^c.
\end{cases}
\end{equation*}
The comparison principle then yields that
$$u(x)\geq\underline{u}(x)\text{~~in~~}\mathbb{R}^n.$$
By the definition of $\rho,$ we have
$\rho(te_n)=d(te_n)$ for any $t\in (0,1),$  and
$$\frac{u(te_n)}{d^s(te_n)}=\frac{u(te_n)}{\rho^s(te_n)}=\beta\frac{u(te_n)}{\underline{u}(te_n)}\geq\beta>0.$$
The proof is complete.
\end{proof}

\section{Regularity}

This section is devoted to  the study of regularity of $(-\triangle)^s_pu.$ We first prove the differentiability of $(-\triangle)_p^su$ under the assumptions of Theorem \ref{differeniability}, then we show that the condition $p>\frac{3}{2-s}$ is optimal by giving a counterexample when $p\leq\frac{3}{2-s}.$

\begin{proof}[\bf Proof of theorem \ref{differeniability}]
For any $x\in \mathbb{R}^n,$ by making  change of variables, $(-\triangle )_p^su(x)$ can be rewritten as
$$(-\triangle)_p^su(x)=\frac{1}{2}C_{n,s,p}\int_{\mathbb{R}^n}\frac{G(u(x)-u(x-y))+G(u(x)-u(x+y))}{|y|^{n+sp}}dy,$$
where $G(t)=|t|^{p-2}t$ for any $t\in\mathbb{R}.$
Note that
\begin{equation}
\begin{aligned}
&\int_{\mathbb{R}^n}\frac{1}{|y|^{n+sp}}\left(\frac{\partial G(u(x)-u(x+y))}{\partial x_i}+ \frac{\partial G(u(x)-u(x-y))}{\partial x_i}\right)dy\\
=&\int_{|y|\leq1}\frac{1}{|y|^{n+sp}}\left(\frac{\partial G(u(x)-u(x+y))}{\partial x_i}+ \frac{\partial G(u(x)-u(x-y))}{\partial x_i}\right)dy\\
&+\int_{|y|> 1} \frac{1}{|y|^{n+sp}}\left(\frac{\partial G(u(x)-u(x+y))}{\partial x_i}+\frac{\partial G(u(x)-u(x-y))}{\partial x_i}\right)dy\\
:=&I_1+I_2.
\end{aligned}
\end{equation}
By a direct calculation, we obtain
\begin{equation}\label{partial-}
\begin{aligned}
\frac{\partial G(u(x)-u(x-y))}{\partial x_i}&=(p-1)|u(x)-u(x-y)|^{p-2}\left(\frac{\partial u(x)}{\partial x_i}
-\frac{\partial u(x-y)}{\partial x_i}\right)\\
&=(p-1)|u(x)-u(x-y)|^{p-2}(v(x)-v(x-y))
\end{aligned}
\end{equation}
and
\begin{equation}\label{partial+}
\frac{\partial G(u(x)-u(x+y))}{\partial x_i}=(p-1)|u(x)-u(x+y)|^{p-2}(v(x)-v(x+y)),
\end{equation}
where $v(x):=\frac{\partial u}{\partial x_i}(x).$

Now, we verify that  $|I_2|<\infty.$ Indeed,
{
\begin{align*}
|I_2|&\leq\int_{|y|> 1} \frac{1}{|y|^{n+sp}}\left(\left|\frac{\partial G(u(x)-u(x+y))}{\partial x_i}\right|+\left|\frac{\partial G(u(x)-u(x-y))}{\partial x_i}\right|\right)dy\\
&\leq C_p\left(\int_{|y|> 1}\frac{|u(x)|^{p-2}|v(x)|+|u(x)|^{p-2}|v(x-y)|}{|y|^{n+sp}}dy \right .\\
&{\phantom{=}\qquad}+\int_{|y|> 1} \frac{|u(x-y)|^{p-2}|v(x)|+|u(x-y)|^{p-2}|v(x-y)|}{|y|^{n+sp}}dy\\
&{\phantom{=}\qquad}+\int_{|y|> 1}\frac{|u(x)|^{p-2}|v(x)|+|u(x)|^{p-2}|v(x+y)|}{|y|^{n+sp}}dy\\
&{\phantom{=}\qquad}\left .+\int_{|y|> 1} \frac{|u(x+y)|^{p-2}|v(x)|+|u(x+y)|^{p-2}|v(x+y)|}{|y|^{n+sp}}dy\right).
\end{align*}
}
It follows from $u\in L_{sp}(\mathbb{R}^n),$  $|\nabla u|\in L_{sp}(\mathbb{R}^n)$ and the  H$\ddot{o}$lder inequality that $|I_2|<\infty.$

For the term $I_1,$  using the Taylor expansion formula, there hold
\begin{equation*}
\begin{aligned}
\frac{\partial G(u(x)-u(x+y))}{\partial x_i}=&(p-1)|\nabla u(x)\cdot y +O(|y|^2)|^{p-2}(-\nabla v(x)\cdot y +O(|y|^2))\\
=&(p-1)|\nabla u(x)\cdot y +O(|y|^2)|^{p-2}(-\nabla v(x)\cdot y)\\
&+(p-1)|\nabla u(x)\cdot y +O(|y|^2)|^{p-2}O(|y|^2)\\
:=&(p-1)J_1+(p-1)J_{2},
\end{aligned}
\end{equation*}
and
\begin{equation*}
\begin{aligned}
\frac{\partial G(u(x)-u(x-y))}{\partial x_i}=&(p-1)|\nabla u(x)\cdot y+O(|y|^2)|^{p-2}(\nabla v(x)\cdot y)\\
&+(p-1)|\nabla u(x)\cdot y +O(|y|^2)|^{p-2}O(|y|^2)\\
:=&(p-1)J_{3}+(p-1)J_4,
\end{aligned}
\end{equation*}
 where the notation $O(|y|^2)$  denotes that there exist some positive constant $C$  such that $|O(|y|^2)|\leq C|y|^2$. Consequently,  there is a  positive constant $C$ such that

\begin{equation}\label{j122}
|J_{2}|+|J_4|\leq C|y|^{p}.
\end{equation}
 Now, we consider the terms $J_1$ and $J_3$.

{\bf Case 1.} $\nabla u(x)\cdot y=0$. Then it follows from the definitions of $J_1$ and $J_3$ that there exist two positive constants $C_1$ and $C_3$ such that
$$|J_1|<C_1|y|^{2p-3}\text{~~and~~} |J_3|<C_3|y|^{2p-3}.$$
{\bf Cases 2.} $\nabla u(x)\cdot y\neq0$. Then we rewrite $J_1$ and $J_3$ respectively as
\begin{equation*}
\begin{aligned}
J_1&=|\nabla u(x)\cdot y +O(|y|^2)|^{p-2}(-\nabla v(x)\cdot y)\\
&=\left(|\nabla u(x)\cdot y +O(|y|^2)|^{p-2}-|\nabla u(x)\cdot y|^{p-2}\right)(-\nabla v(x)\cdot y)
-|\nabla u(x)\cdot y|^{p-2}\nabla v(x)\cdot y,\\
\end{aligned}
\end{equation*}
and
\begin{equation*}
J_3=\left(|\nabla u(x)\cdot y +O(|y|^2)|^{p-2}-|\nabla u(x)\cdot y|^{p-2}\right)(\nabla v(x)\cdot y)
+|\nabla u(x)\cdot y|^{p-2}\nabla v(x)\cdot y.
\end{equation*}
It follows  that
\begin{equation*}
\begin{aligned}
|J_1+J_3|\leq& C\left(|\nabla u(x)\cdot y +O(|y|^2)|^{p-2}-|\nabla u(x)\cdot y|^{p-2}\right)|\nabla v(x)\cdot y|\\
\leq&C|\nabla u(x)\cdot y|^{p-4}|2O(|y|^2)\nabla u(x)\cdot y +O(|y|^4)||\nabla v(x)\cdot y|\\
\leq&C|y|^p.
\end{aligned}
\end{equation*}
To summary, we conclude that there exists a positive constant $C$ independent of $y$ such that
$$|J_1+J_3|+|J_2|+|J_4|\leq C(|y|^{2p-3}+|y|^p)\text{~~for any~} y\in B_1(0).$$
The assumption $p>\frac{3}{2-s}$ further implies
$$|I_1|\leq(p-1)\int_{|y|\leq1}\frac{1}{|y|^{n+sp}}(|J_{2}|+|J_{4}|+|J_1+J_3|)<\infty,$$
that is,
$$\int_{\mathbb{R}^n}\frac{1}{|y|^{n+sp}}\left|\frac{\partial G(u(x)-u(x+y))}{\partial x_i}+ \frac{\partial G(u(x)-u(x-y))}{\partial x_i}\right|dy<\infty \text{~~for any~}x\in \mathbb{R}^n.$$
By exchanging the order of integration and differentiation, we derive that $(-\triangle )_p^su$ is differentiable in  $\mathbb{R}^n,$ and then we conclude
$(-\triangle )_p^s u \in C(\mathbb{R}^n)$ by exchanging  the order of integration  and limit. The proof is complete.
\end{proof}

Theorem \ref{differeniability} verifies that in the case $p>2,$  if one assumes in addition that $p>\frac{3}{2-s},$
 then $(-\triangle )_p^su \in C^1(\mathbb{R}^n)$ for any $u$ satisfying $u\in C_{\text{loc}}^3(\mathbb{R}^n)\cap L_{sp}(\mathbb{R}^n)$ and  $|\nabla u|\in L_{sp}(\mathbb{R}^n).$ It seems from the proof of Theorem \ref{differeniability} that $p>\frac{3}{2-s}$ is a technical assumption. While, the counterexample in Theorem \ref{opposite} shows that this  condition is optimal to ensure $(-\triangle )_p^s u \in C(\mathbb{R}^n)$ for any $u$ satisfying $u\in C_{\text{loc}}^3(\mathbb{R}^n)\cap L_{sp}(\mathbb{R}^n)$ and  $|\nabla u|\in L_{sp}(\mathbb{R}^n)$.

\begin{proof}[\bf Proof of Theorem \ref{opposite}]
By virtue of the definition, we have
$$(-\triangle)_p^su(x)=\frac{1}{2}C_{s,p}\int_{-\infty}^{+\infty}\frac{G(u(x)-u(x+y))+G(u(x)-u(x-y))}{|y|^{1+sp}}dy.$$
For the convenience of writing, we set
\begin{equation}\label{Fxy}
F(x,y):=\frac{G(u(x)-u(x+y))+G(u(x)-u(x-y))}{|y|^{1+sp}}.
\end{equation}
It follows from a straightforward calculation that
{
\begin{align}\label{Fx}
\frac{\partial F(x,y)}{\partial x}
=&\frac{(p-1)}{|y|^{1+sp}}\big(|u(x)-u(x+y)|^{p-2}(u^\prime(x)-u^\prime(x+y))\notag\\
&+|u(x)-u(x-y)|^{p-2}(u^\prime(x)-u^\prime(x-y))\big)\notag\\
=&\frac{(p-1)}{|y|^{1+sp}}|\eta(x)x^2-\eta(x+y)(x+y)^2|^{p-2}\notag\\
&\times\big[\eta^\prime(x)x^2+2\eta(x)x-\eta^\prime(x+y)(x+y)^2-2\eta(x+y)(x+y)\big]\notag\\
&+\frac{(p-1)}{|y|^{1+sp}}|\eta(x)x^2-\eta(x-y)(x-y)^2|^{p-2}\notag\\
&\times\big[\eta^\prime(x)x^2+2\eta(x)x-\eta^\prime(x-y)(x-y)^2-2\eta(x-y)(x-y)\big].
\end{align}
}
Let
$$f(x,y)=|\eta(x)x^2-\eta(x+y)(x+y)^2|^{p-2}\big[\eta^\prime(x)x^2+2\eta(x)x-\eta^\prime(x+y)(x+y)^2-2\eta(x+y)(x+y)\big],$$
then we can rewrite \eqref{Fx} as
$$\frac{\partial F(x,y)}{\partial x}=\frac{(p-1)}{|y|^{1+sp}}[f(x,y)+f(x,-y)].$$
Note that for any $0<x<\frac{1}{8},$  there hold
{
\begin{equation}\label{example2}
\begin{aligned}
&\int_{-\infty}^{+\infty}\frac{\partial F(x,y)}{\partial x}dy\\
&=(p-1)\int_{-\infty}^{+\infty}\frac{f(x,y)+f(x,-y)}{|y|^{1+sp}}dy\\
&=2(p-1)\int_{0}^{+\infty}\frac{f(x,y)+f(x,-y)}{|y|^{1+sp}}dy\\
&=2(p-1)\left(\int_{0}^{\frac{1}{2}}\frac{f(x,y)+f(x,-y)}{|y|^{1+sp}}dy+\int^{\frac{5}{2}}_{\frac{1}{2}}\frac{f(x,y)+f(x,-y)}{|y|^{1+sp}}dy\right .\\
&{\phantom{=}\qquad\qquad\quad}\left .+\int_{\frac{5}{2}}^{\infty}\frac{f(x,y)+f(x,-y)}{|y|^{1+sp}}dy\right)\\
&:=2(p-1)(I_1+I_2+I_3).
\end{aligned}
\end{equation}
}
For $I_3,$ in view of $y>\frac{5}{2}$ and $0<x<\frac{1}{8},$ there hold $|x-y|>2$ and $|x+y|>2,$ which along with the properties of $\eta$ implies that
$$f(x,y)+f(x,-y)=4x^{2p-3}.$$
Hence, we have
\begin{equation}\label{example4}
I_3=\int_{\frac{5}{2}}^\infty\frac{4x^{2p-3}}{y^{1+sp}}dy= C_1x^{2p-3},
\end{equation}
where $C_1$ is a positive constant independent of $x.$

For $I_2,$ thanks to $\frac{1}{2}<y<\frac{5}{2}$ and $0<x<\frac{1}{8}$, there exists a positive constant $C_2$ independent of $x$ such that
\begin{equation}\label{example5}
|I_2|<\infty.
\end{equation}

{It remains to} estimate the term $I_1.$ By virtue of $0<x<\frac{1}{8}$ and $0<y<\frac{1}{2},$ we see $|x-y|<1$ and
$|x+y|<1,$ which together with the properties of $\eta$ yields
$$f(x,y)+f(x,-y)=2y|y^2-2xy|^{p-2}-2y|y^2+2xy|^{p-2}.$$
Therefore, for any $x\in (0,\frac{1}{8}),$ there hold
\begin{align}\label{example1}
I_1&=\int_0^{\frac{1}{2}}\frac{2y|y^2-2xy|^{p-2}-2y|y^2+2xy|^{p-2}}{|y|^{1+sp}}dy\notag\\
&=\int_0^{\frac{1}{2}}\frac{2}{y^{2+sp-p}}(|y-2x|^{p-2}-|y+2x|^{p-2})dy\notag\\
&=2\int_{0}^{\frac{1}{2x}}\frac{x}{(xz)^{2+sp-p}}(|xz-2x|^{p-2}-|xz+2x|^{p-2})dz\notag\\
&=\frac{2}{x^{sp-2p+3}}\int_{0}^{\frac{1}{2x}}\frac{|z-2|^{p-2}-|z+2|^{p-2}}{z^{2+sp-p}}dz\notag\\
&=\frac{2}{x^{sp-2p+3}}\int_{0}^2\frac{(2-z)^{p-2}-(z+2)^{p-2}}{z^{2+sp-p}}dz\notag\\
&{\phantom{=}}+\frac{2}{x^{sp-2p+3}}\int_{2}^{\frac{1}{2x}}\frac{(z-2)^{p-2}-(z+2)^{p-2}}{z^{2+sp-p}}dz\notag\\
&=\frac{2}{x^{sp-2p+3}}\int_{0}^2\frac{(2-z)^{p-2}-2^{p-2}+2^{p-2}-(z+2)^{p-2}}{z^{2+sp-p}}dz\notag\\
&{\phantom{=}}+\frac{2}{x^{sp-2p+3}}\int_{2}^{\frac{1}{2x}}\frac{(z-2)^{p-2}-(z+2)^{p-2}}{z^{2+sp-p}}dz\notag\\
&=\frac{2}{x^{sp-2p+3}}\int_{0}^2\frac{2(2-p)2^{p-3}z+O(|z|^2)}{z^{2+sp-p}}dz\notag\\
&{\phantom{=}}+\frac{2}{x^{sp-2p+3}}\int_{2}^{\frac{1}{2x}}\frac{(z-2)^{p-2}-(z+2)^{p-2}}{z^{2+sp-p}}dz\notag\\
&<\infty,
\end{align}
which along with  \eqref{example2},\ \eqref{example4} and \eqref{example5}  shows that for any fixed $x\in (0,\frac{1}{8}),$
there holds
$$\left|\int_{-\infty}^{+\infty}\frac{\partial F(x,y)}{\partial x}dy\right|<\infty.$$
By exchanging the order of integration and differentiation, we derive that $((-\triangle)_p^su)^\prime$ is well-defined for any
$x\in (0,\frac{1}{8}),$ and
 $$((-\triangle)_p^su)^\prime(x)=
 \frac{1}{2}C_{s,p}\int_{-\infty}^{+\infty}\frac{\partial F(x,y)}{\partial x}dy.$$
 If $p<\frac{3}{2-s},$ then \eqref{example1} implies that
 \begin{align*}
 I_1&=\frac{2}{x^{sp-2p+3}}\int_{0}^2\frac{(2-z)^{p-2}-(z+2)^{p-2}}{z^{2+sp-p}}dz\\
&{\phantom{=}}+\frac{2}{x^{sp-2p+3}}\int_{2}^{\frac{1}{2x}}\frac{(z-2)^{p-2}-(z+2)^{p-2}}{z^{2+sp-p}}dz\\
&\leq\frac{2}{x^{sp-2p+3}}\int_{0}^2\frac{(2-z)^{p-2}-2^{p-2}}{z^{2+sp-p}}dz\\
&=\frac{2}{x^{sp-2p+3}}\int_{0}^2\frac{(2-p)2^{p-3}z+O(z^2)}{z^{2+sp-p}}dz\\
&\leq -Cx^{2p-sp-3},
 \end{align*}
 which verifies that
$$\lim_{x\to0+}I_1=-\infty.$$
Consequently, we conclude that
\begin{align*}
\lim_{x\to 0+}((-\triangle)_p^su)^\prime(x)&=\lim_{x\to 0+}\frac{1}{2}C_{s,p}\int_{-\infty}^{+\infty}\frac{\partial F(x,y)}{\partial x}dy\\
&=\lim_{x\to 0+}C_{s,p}(p-1)(I_1+I_2+I_3)\\
&=-\infty,
\end{align*}
that is, \eqref{example} holds.

In the case $p=\frac{3}{2-s},$ we first prove that $\left((-\triangle)_p^su\right)^\prime(0)=0.$ In fact,
{
\begin{align*}
\left((-\triangle)^s_pu\right)^\prime_+(0)&=\lim_{x\to0+}\frac{(-\triangle)^s_pu(x)-(-\triangle)^s_pu(0)}{x}\\
&=\frac{1}{2}C_{s,p}\lim_{x\to0+}\frac{1}{x}\int_{-\infty}^{+\infty}F(x,y)-F(0,y)dy\\
&=C_{s,p}\lim_{x\to0+}\frac{1}{x}\int_{0}^{+\infty}F(x,y)-F(0,y)dy\\
&=C_{s,p}\Big(\lim_{x\to0+}\int_{0}^{\frac{1}{2}}\frac{F(x,y)-F(0,y)}{x}dy\\
&{\phantom{=}\qquad\qquad}+\lim_{x\to0+}\int_{\frac{1}{2}}^{\frac{5}{2}}\frac{F(x,y)-F(0,y)}{x}dy\\
&{\phantom{=}\qquad\qquad}+\lim_{x\to0+}\int_{\frac{5}{2}}^{+\infty}\frac{F(x,y)-F(0,y)}{x}dy\Big)\\
&:=C_{s,p}(J_1+J_2+J_3)
\end{align*}
}
For $J_3,$ in view of  $y>\frac{5}{2}$ and $0<x<\frac{1}{8},$ there hold $|x-y|>2$ and $|x+y|>2,$ which along with the properties of $\eta$ and $sp=2p-3$ implies that
$$F(x,y)=\frac{2x^{2p-2}}{y^{2p-2}}~~ \text{and}~F(0,y)=0.$$
It then follows that
\begin{equation}
J_3=2\lim_{x\to0+}x^{2p-3}\int_{\frac{5}{2}}^{+\infty}\frac{1}{y^{2p-2}}dy=0.
\end{equation}
For $J_2,$ by exchanging the order of integration and limit, we have
{
\begin{equation}\label{JJ}
\begin{aligned}
J_2&=\int_{\frac{1}{2}}^{\frac{5}{2}}\frac{\partial F}{\partial x}(0,y)dy\\
&=(p-1)\int_{\frac{1}{2}}^{\frac{5}{2}}\left(\frac{\left|\eta(y)y^2\right|^{p-2}(-y^2\eta^\prime(y)-2y\eta(y))}{y^{2p-2}}\right .\\
&{\phantom{=}}\left.
~~~~~~~~~~~~~~+\frac{\left|\eta(-y)y^2\right|^{p-2}(-y^2\eta^\prime(-y)+2y\eta(-y))}{y^{2p-2}}\right)dy
\end{aligned}
\end{equation}
}
Since $\eta(y)=\eta(-y)$ in $\mathbb{R},$ there holds $\eta^\prime(y)=-\eta^\prime(-y),$ which along with \eqref{JJ} implies
$$J_2=0.$$
As for $J_1$, we see
\begin{align*}
J_1&=\lim_{x\to0+}\frac{1}{x}\int_{0}^{\frac{1}{2}}F(x,y)-F(0,y)dy\\
&=\lim_{x\to0+}\frac{1}{x}\int_{0}^{\frac{1}{2}}\frac{\partial F}{\partial x}(0,y)x+O(x^2)dy\\
&=\int_{0}^{\frac{1}{2}}\frac{\partial F}{\partial x}(0,y)dy\\
&=(p-1)\int_{0}^{\frac{1}{2}}\frac{|y^2|^{p-2}(-2y)+|y^2|^{p-2}(2y)}{y^{2p-2}}dy\\
&=0.
\end{align*}
To summary, we conclude that
 $$\left((-\triangle)_p^su\right)_+^\prime(0)=0.$$
 Similarly, we can prove
 $$\left((-\triangle)_p^su\right)_-^\prime(0)=0.$$
It then follows that
$$\left((-\triangle)_p^su\right)^\prime(0)=0.$$

On the other hand,
\eqref{example4} implies that
\begin{equation}\label{ij3}
\lim_{x\to 0+}I_3=0.
\end{equation}
Thanks to \eqref{example5}, by exchanging the order of integration and limit, we obtain
\begin{equation}\label{ij2}
\lim_{x\to0+}I_2=0.
\end{equation}
A similar calculation to \eqref{example1} implies that for any $x\in (0,\frac{1}{8}),$
\begin{equation}\label{I1}
\begin{aligned}
I_1&=2\int_{0}^{\frac{1}{2x}}\frac{|z-2|^{p-2}-|z+2|^{p-2}}{z^{p-1}}dz\\
&=2\int_{0}^2\frac{(2-z)^{p-2}-(z+2)^{p-2}}{z^{p-1}}dz\\
&\quad+2\int_{2}^{\frac{1}{2x}}\frac{(z-2)^{p-2}-(z+2)^{p-2}}{z^{p-1}}dz\\
&<2\int_{0}^2\frac{(2-z)^{p-2}-(z+2)^{p-2}}{z^{p-1}}dz\\
&<2\int_{1}^2\frac{(2-z)^{p-2}-(z+2)^{p-2}}{z^{p-1}}dz\\
&\leq \frac{2(1-3^{p-2})}{2^{p-1}}\\
&<0.
\end{aligned}
\end{equation}
To summary, we derive that
\begin{align*}
\lim_{x\to0+}((-\triangle)_p^su)^\prime(x)&=
 \frac{1}{2}C_{s,p}\lim_{x\to0+}\int_{-\infty}^{+\infty}\frac{\partial F(x,y)}{\partial x}dy\\
 &=C_{s,p}(p-1)\lim_{x\to0+}(I_1+I_2+I_3)\\
 &<C_{s,p}(p-1)\frac{2(1-3^{p-2})}{2^{p-1}}\\
 &<0.
\end{align*}
The proof is complete.
\end{proof}


\begin{appendices}
\section{}
In this Appendix, we list some results in \cite{iannizzotto} that were used in the proof of Theorem \ref{hopf}.
The first one is the  weak comparison principle.

\begin{lemma}
Let $\Omega\subset\mathbb{R}^n$ be a bounded domain. Assume $u,v\in \widetilde{W}^{s,p}(\Omega)$ satisfy, in the weak sense,
\begin{equation*}
\begin{cases}
(-\triangle)_p^su\geq(-\triangle)_p^sv &\text{~in~}\Omega,\\
u\geq v&\text{~in}~\mathbb{R}^n\backslash\Omega.
\end{cases}
\end{equation*}
Then $$u\geq v\quad\text{a.e.~in}~\Omega.$$
\end{lemma}

Another key ingredient is the following ``change of variables" lemma.

\begin{lemma}\label{changeofvariable}
Let $\Psi$ be a $C^{1,1}$ diffeomorphism of $\mathbb{R}^n$ such that $\Psi=I$ in $B_r^c(0),$ $r>0.$
Then the function $v(x)=(\Psi^{-1}(x)\cdot e_n)_+^s$ belongs to $\widetilde{W}_{\text{loc}}^{s,p}(\mathbb{R}^n)$
and is a weak solution of
$$(-\triangle)_p^sv=f \quad \text{in}~\Psi(\mathbb{R}^n_+),$$
with
$$||f||_\infty\leq C\left(||D\Psi||_\infty,||D\Psi^{-1}||_\infty,r\right)||D^2\Psi||_\infty,$$
where $C\left(||D\Psi||_\infty,||D\Psi^{-1}||_\infty,r\right)$ is a positive constant.
Moreover,
\begin{equation}\label{change}
\lim_{\epsilon\to0}C_{n,s,p}\int_{\{|\Psi^{-1}(x)-\Psi^{-1}(y)|>\epsilon\}}
\frac{|v(x)-v(y)|^{p-2}(v(x)-v(y))}{|x-y|^{n+sp}}dy=f\text{~in~} L^1_{\text{loc}}(\Psi(\mathbb{R}^n_+)).
\end{equation}
\end{lemma}
\begin{remark}
The equality \eqref{change} follows from the proof of ``change of variables" lemma.
\end{remark}

The following lemma implies that the point-wise solution is also a weak solution.

\begin{lemma}\label{generaldef}
Let $u\in \widetilde{W}^{s,p}_{\text{loc}}(\Omega)$ and $D$ denote the diagonal of $\mathbb{R}^n\times\mathbb{R}^n.$ For any $\epsilon>0,$ assume $A_\epsilon\subset\mathbb{R}^n\times\mathbb{R}^n$ is a neighborhood of $D$ and satisfies
\begin{itemize}
\item[(i)] $(x,y)\in A_\epsilon$ for all $(y,x)\in A_\epsilon,$
\item[(ii)]$\mathop{\sup}\limits_{x\in A_\epsilon}\text{dist}(x,D)\to 0 $ as $\epsilon\to0$.
\end{itemize}
For any $x\in\mathbb{R}^n,$ we set $A_\epsilon(x)=\{y\in\mathbb{R}^n|(x,y)\in A_\epsilon\}$ and
$$g_\epsilon(x)=C_{n,s,p}\int_{A^c_\epsilon(x)}\frac{|u(x)-u(y)|^{p-2}(u(x)-u(y))}{|x-y|^{n+sp}}dy.$$
If $g_\epsilon\to f$ in $L^1_{\text{loc}}(\Omega),$ then $u$ is a weak solution of
$$(-\triangle)^s_pu=f\text{~~in~}\Omega.$$
\end{lemma}

\end{appendices}



\begin{thebibliography}{99}

\bibitem{ab}
G. Alberti,  G. Bellettini, A nonlocal anisotropicmodel for phase transitions I: the optimal profile problem.
Math. Ann. 310 (1998), 527-560.


\bibitem{BCF1}
C. Bjorland,  L. Caffarelli,  A. Figalli,  Non-local gradient dependent operators. Adv. Math. 230 (2012),  1859-1894.
\bibitem{BCF2}
 C. Bjorland,  L. Caffarelli,  A. Figalli,  Nonlocal tug-of-war and the infinity fractional Laplacian. Comm. Pure Appl. Math. 65 (2012),  337-380.

\bibitem{brasco}
 L. Brasco, E. Lindgren, Armin Schikorra, Higher H$\ddot{o}$der regularity for the fractional p-Laplacian in the superquadratic case, arXiv:1711.09835.
 \bibitem{bcps}
  C. Brandle,  E. Colorado, A. de Pablo,  U. Sanchez,  A concave-convex elliptic problem involving the fractional Laplacian. Proc. Roy. Soc. Edinburgh Sect. A 143 (2013),  39-71.

  \bibitem{bgr}
K. Bogdan, T. Grzywny,  M. Ryznar, Heat kernel estimates for the
fractional Laplacian with Dirichlet conditions. Ann. of Prob. 38 (2010),
1901-1923.
\bibitem{cfy}
 W. Chen, Y. Fang, R. Yang, Liouville theorems involving the fractional Laplacian
on a half space, Advances in Math. 274 (2015), 167-198.
 \bibitem{chenli1}
 W. Chen, C. Li, Maximum principles for the fractional p-Laplacian and symmetry of solutions, arXiv:1705.04891.
 \bibitem{chenli2}
 W. Chen, C. Li, A Hopf type lemma for fractional equations, arXiv:1705.04889.

 \bibitem{chenlizhang}
 W. Chen,  Y. Li,  R. Zhang,  A direct method of moving spheres on fractional order equations. J. Funct. Anal. 272 (2017), 4131-4157.

\bibitem{chenlili}
W. Chen, C. Li,  G. Li, Maximum principles for a fully nonlinear fractional order equation and symmetry of solutions. Calc. Var. Partial Differential Equations 56 (2017), no. 2, Art. 29, 18 pp.

\bibitem{chenliliyan1}
 W. Chen,  C. Li,  Y. Li,  A direct method of moving planes for the fractional Laplacian. Adv. Math. 308 (2017), 404-437.

\bibitem{chenliliyan2}
 W. Chen,  C. Li,  Y. Li, A direct blowing-up and rescaling argument on nonlocal elliptic equations. Internat. J. Math. 27 (2016),  1650064, 20 pp.

 \bibitem{chenlima}
 W. Chen, Y. Li,  P. Ma, The Fractional Laplacian, a book to be published by
World Scientific Publishing Co. 2017.


 \bibitem{chenzhu}
  W. Chen, J. Zhu,  Indefinite fractional elliptic problem and Liouville theorems. J. Differential Equations 260 (2016),  4758-4785.

  \bibitem{cs}
  L. Caffarelli, L. Silvestre, An extension problem related to the fractional Laplacian, Comm. PDE. 32(2007), 1245-1260.
\bibitem{ck}
 A. Di Castro,  T. Kuusi,  G. Palatucci, Local behavior of fractional p-minimizers. Ann. Inst. H. Poincare Anal. Non Lineaire 33 (2016),  1279-1299.
\bibitem{dpv}
E. Di Nezza, G. Palatucci,  E. Valdinoci,  Hitchhiker's guide to the fractional Sobolev spaces. Bull. Sci. Math. 136 (2012),  521-573.
\bibitem{fall-jarohs}
 M. M. Fall, S. Jarohs,  Overdetermined problems with fractional Laplacian. ESAIM Control Optim. Calc. Var. 21 (2015), no. 4, 924-938.

\bibitem{grecoservadei}
 A. Greco, R. Servadei,  Hopf's lemma and constrained radial symmetry for the fractional Laplacian. Math. Res. Lett. 23 (2016), 863-885.
\bibitem{iannizzotto}
 A. Iannizzotto,  S. Mosconi,  M. Squassina,  Global H$\ddot{o}$der regularity for the fractional p-Laplacian. Rev. Mat. Iberoam. 32 (2016),  1353-1392.
 \bibitem{ishiinakamura}
 H. Ishii, G. Nakamura, A class of integral equations
and approximation of p-Laplace equations, Calc. Var. Partial Differential
Equations 37 (2010),  485-522.
\bibitem{jinli}
L. Jin, Y. Li, A Hopf's Lemma and the Boundary Regularity for the Fractional P-Laplacian, arXiv:1711.02707.
\bibitem{ms}
 S. Mosconi,  M. Squassina,  Recent progresses in the theory of nonlinear nonlocal problems. Bruno Pini Mathematical Analysis Seminar 2016, 147-164, Bruno Pini Math. Anal. Semin., 2016, Univ. Bologna, Alma Mater Stud., Bologna, 2016.
\bibitem{otonserra}
X. Ros-Oton, J. Serra, The Dirichlet problem for the fractional
Laplacian: regularity up to the boundary. J. Math. Pures Appl. 9
(2014), 275-302.

\bibitem{silvestre}
 L. Silvestre,  Regularity of the obstacle problem for a fractional power of the Laplace operator. Comm. Pure Appl. Math. 60 (2007),  67-112.
\bibitem{sirevaldinoci}
Y. Sire, E. Valdinoci,  Fractional Laplacian phase transitions and boundary reactions: a geometric inequality
and a symmetry result. J. Funct. Anal. 256 (2009), 1842-1864.

\bibitem{valdinoci}
E. Valdinoci,  From the long jump random walk to the fractional Laplacian. Bol. Soc. Esp. Mat. Apl. Se
MA 49 (2009), 33-44.
\bibitem{zccy}
 R. Zhuo, W. Chen, X. Cui,  Z. Yuan, A Liouville theorem for the fractional
Laplacian,  arXiv:1401.7402.
\end{thebibliography}
\end{document}